\documentclass[review]{elsarticle}
\usepackage{amsmath}
\usepackage{amsfonts}
\usepackage{amssymb}
\usepackage{amsthm}

\usepackage{lineno,hyperref}
\modulolinenumbers[5]

\journal{Journal of \LaTeX\ Templates}

\newtheorem{theorem}{\sc\bf Theorem}[section]
\newtheorem{example}[theorem]{\sc\bf Example}
\newtheorem{lemma}[theorem]{\sc\bf Lemma}
\newtheorem{corollary}[theorem]{\sc\bf Corollary}
\newtheorem{remark}[theorem]{\sc\bf Remark}
\newtheorem{definition}[theorem]{\sc\bf Definition}

\renewcommand{\d}[1]{\ensuremath{\operatorname{d}\!{#1}}}

\begin{document}

\begin{frontmatter}

\title{On the Ulam Type Stability of Nonlinear Volterra Integral Equations}


\author[mymainaddress]{Süleyman Öğrekçi\corref{mycorrespondingauthor}}
\address[mymainaddress]{Science \& Arts Faculty, Departmet of Mathematics, Amasya University, Turkey}
\cortext[mycorrespondingauthor]{Corresponding author}
\ead{suleyman.ogrekci@amasya.edu.tr}

\author[mysecondaryaddress]{Yasemin Başcı}
\address[mysecondaryaddress]{Science \& Arts Faculty, Departmet of Mathematics, Bolu Abant Izzet Baysal University, Turkey}
\ead{basci\_y@ibu.edu.tr}

\author[mythirdaddress]{Adil Mısır}
\address[mythirdaddress]{Faculty of Science, Departmet of Mathematics, Gazi University, Turkey}
\ead{adilm@gazi.edu.tr}

\begin{abstract}
In this paper, we examine the Hyers-Ulam and Hyers-Ulam-Rassias stability of solutions of a general class of nonlinear Volterra integral equations. By using a fixed point alternative and improving a technique commonly used in similar problems, we extend and improve some well-known results on this problem. We also provide some examples visualizing the improvement of the results mentioned.
\end{abstract}


\end{frontmatter}


\section{Introduction}

Hyers–Ulam stability, initiated with a speech of S. M. Ulam \cite{book-ulam} at Wisconsin University, is a concept that provides an approximate solution for the exact solution in a simple form for differential equations. Ulam posed the following problem: ``Under what conditions does there exists an homomorphism near an approximately homomorphism of a complete metric group?'' More precisely: Given a metric group $\left(G, \cdot, d \right)$, a number $\varepsilon>0$ and a mapping $f:G\rightarrow G$ satisfying the inequality $$d\left( f(xy), f(x)f(y) \right)<\varepsilon$$ for all $x,y\in G$, does there exist a homomorphism $g$ of $G$ and a constant $K$, depending only on $G$, such that $$d\left( f(x), g(x) \right)\leq K\varepsilon$$ for all $x\in G$? In the presence of affirmative answer, the equation $g(xy)=g(x)g(y)$ of the homomorphism is called \emph{stable}, see \cite{book-ulam} for details. One year later, Hyers \cite{hyers} gave an answer to this problem for linear functional equations on Banach spaces and showed that the additive functional equation $f(x+y)=f(x)+f(y)$ is stable in the sense of Ulam with the following celebrated result.
\begin{theorem}[\cite{hyers}]
	Let $E_1, E_2$ be real Banach spaces and $\varepsilon>0$. Then, for each mapping $f:E_1\rightarrow E_2$ satisfying $$\left\Vert f(x+y)-f(x)-f(y) \right\Vert\leq\varepsilon$$ for all $x,y\in E_1$, the limit $$L(x):=\lim\limits_{n\to\infty}2^{-n}f\left(2^nx\right)$$ exists for all $x\in E_1$ and $L:E_1\to E_2$ is the unique additive mapping satisfying $$\left\|f(x)-L(x)\right\|\leq\varepsilon$$ for all $x\in E_1$.
\end{theorem}

After this result of Hyers, a new concept of stability for functional equations established, called today Hyers-Ulam stability, and many papers devoted to this subject (see for example \cite{cadariu, forti, brzdek, petru, b1,b2,b3}). In 1978, T. Rassias \cite{rassias} provided a remarkable generalization with the following well-known result, which known as Hyers-Ulam-Rassias stability today.
\begin{theorem}[\cite{rassias}]
	Let $E_1, E_2$ be Banach spaces, $\theta\in(0,\infty)$ and $p\in[0,1)$. Then, for each mapping $f:E_1\rightarrow E_2$ satisfying $$\left\Vert f(x+y)-f(x)-f(y) \right\Vert\leq\varepsilon$$ for all $x,y\in E_1$, there exists a unique additive mapping $L:E_1\to E_2$ satisfying $$\left\|f(x)-L(x)\right\|\leq\frac{2\theta}{2-2^p}\left\|x\right\|^p$$ for all $x\in E_1$. Moreover, if $f(tx)$ is continuous for all $t\in\mathbb{R}$ and each fixed $x\in E_1$, then $L$ is $\mathbb{R}$-linear.
\end{theorem}

Apart from functional equations, this concept of stability has applied to various kind of equations such as differential equations, integral equations and integrodifferential equations. Stability problem of differential equations in the sense of Hyers-Ulam was initiated by the papers of Obloza \cite{obloza1, obloza2}. Later Alsina and Ger \cite{alsina}, Takahasi et. al. \cite{takahasi1,takahasi2, takahasi3} provided remarkable results in this topic. After these pioneering works, a large number of papers devoted to this subject have been published (see for example \cite{jung,ben1, ben2,ben3,ben4,Popa2016, Shen2017,deOliveira2018} and references therein).

Volterra integral equations have been studied extensively since the four fundamental papers of Vito Volterra in 1896, and specially since 1913 when Volterra's book ``Leçons sur les \'Equations Int\'egrales et les \'Equations Int\'egro-diff\'erentielles'' appeared. For a continuous function $f$ and a fixed constant $a\in\mathbb{R}$, the integral equation
\begin{equation}\label{eq1}
y(t)=\int_{a}^{t}f(s,y(s))\d{s}
\end{equation}
is called Volterra integral equation of second kind. Jung \cite{jung2} studied Hyers-Ulam and Hyers-Ulam-Rassias stability of the equation \eqref{eq1} by using a fixed point alternative. Further, Castro \cite{castro} investigated Ulam type stability criteria for  the equations
\begin{equation}\label{eq2}
y(t)=\int_{a}^{t}f(t,s,y(s))\d{s}.
\end{equation}
Later Gachpazan and Baghani \cite{gachpazan} and Akkouchi \cite{akkouchi} considered the stability of Volterra integral equations $$y(t)=f(t)+\lambda\int_{t_0}^{t}k(t,s)y(s)\d{s}\quad\textrm{and}\quad y(t)=f(t)+\lambda\int_{t_0}^{t}f(t,y,f(y))\d{s}$$ successively in the sense of Ulam. There are also some works in the literature on the stability problem of Volterra integral equation with delay, for example Castro \cite{castro2} studied this problem for the equation $$y(t)=f(t)+\Psi\left(\int_{t_0}^{t}f(t,s, y(s),y(\alpha(s))\d{s}\right).$$

Now, we give explicit definitions of the Hyers-Ulam and Hyers-Ulam-Rassias stability concepts for Volterra integral equation \eqref{eq1}, these definitions can be adapted to the equations \eqref{eq2} and other equations mentioned above in similar ways.
\begin{definition}\label{def:ulam}
	If for each function $y(t)$ satisfying
	\begin{equation*}
		\left| y(t)-\int_{a}^{t}f(s,y(s))\d{s} \right|\leq\varepsilon,
	\end{equation*}
	for all $t$ and some $\varepsilon\geq0$, there exists a solution $y_0(t)$ of \eqref{eq1} and a constant $K>0$ independent of $y$ and $y_0$ satisfying
	\begin{equation*}
		\left| y(t)-y_0(t) \right|\leq K\varepsilon
	\end{equation*}
	for all $t$, then the Volterra integral equation \eqref{eq1} is said to be stable in the sense of Hyers-Ulam \cite{jung2}.
\end{definition}

\begin{definition}\label{def:rassias}
	If the statement of Definition \ref{def:ulam} is true after replacing the constant $\varepsilon$ with the function $\varphi(t)\geq0$, where this function does not depend on $y$ and $y_0$, then the Volterra integral equation \eqref{eq1} is said to be stable in the sense of Hyers-Ulam-Rassias \cite{jung2}.
\end{definition}

In \cite{jung2}, Jung proved following remarkable results on the Ulam type stability of the solutions of the Volterra integral equation \eqref{eq1}.

\begin{theorem}[\cite{jung2}]
	Let the function $f:[a,b]\times\mathbb{C}\to\mathbb{C}$ be a continuous function satisfying Lipschitz condition with the Lipschitz constant $L$ on $[a-r,a+r]$ with
	\begin{equation}\label{eq:jung1}
	Lr<1,
	\end{equation}
	then the Volterra integral equation \eqref{eq1} has Hyers-Ulam stability on $[a-r,a+r]$.
\end{theorem}

\begin{theorem}[\cite{jung2}]
	Let the function $f:[a,b]\times\mathbb{C}\to\mathbb{C}$ be a continuous function satisfying Lipschitz condition with the Lipschitz constant $L$ on $[a,b]$. If the function $\varphi:[a,b]\to(0,\infty)$ as in Definition \ref{def:rassias} satisfies
	\begin{equation}\label{eq:jung3}
	\left| \int_{a}^{t}\varphi(s)\d{s} \right|\leq K\varphi(t)
	\end{equation}
	and
	\begin{equation}\label{eq:jung2}
	KL<1
	\end{equation}
	for a positive constant $K$, then the Volterra integral equation \eqref{eq1} has Hyers-Ulam-Rassias stability on $[a,b]$.
\end{theorem}

In \cite{castro}, Castro and Ramos obtained the following similar results on the stability problem of the Volterra integral equation \eqref{eq2}.

\begin{theorem}[\cite{castro}]\label{castro1}
	Let the function $f:[a,b]\times[a,b]\times\mathbb{C}\to\mathbb{C}$ be a continuous function satisfying Lipschitz condition with the Lipschitz constant $L$ on $[a-r,a+r]$ satisfying \eqref{eq:jung1},
	then the Volterra integral equation \eqref{eq2} has Hyers-Ulam stability on $[a-r,a+r]$.
\end{theorem}

\begin{theorem}[\cite{castro}]\label{castro2}
	Let the function $f:[a,b]\times[a,b]\times\mathbb{C}\to\mathbb{C}$ be a continuous function satisfying Lipschitz condition with the Lipschitz constant $L$ on $[a,b]$. If the function $\varphi:[a,b]\to(0,\infty)$ as in Definition \ref{def:rassias} satisfies \eqref{eq:jung3} and \eqref{eq:jung2} for a positive constant $K$, then the Volterra integral equation \eqref{eq2} has Hyers-Ulam-Rassias stability on $[a,b]$.
\end{theorem}

We remark that all the results mentioned above heavily depend on the conditions \eqref{eq:jung1}, \eqref{eq:jung3} and \eqref{eq:jung2}. The main purpose of this paper is to examine the Hyers-Ulam and Hyers-Ulam-Rassias stability of Volterra integral equations given by \eqref{eq1} and \eqref{eq2}. We will improve a technique seen in the work of C\u{a}dariu and Radu \cite{cadariu2} and also used in the papers \cite{jung2} and \cite{castro}, this improvement will enable us to obtain more accurate results for the problem considered and to improve the results given in \cite{jung2} and \cite{castro}. More precisely, we will show that the conditions \eqref{eq:jung1}, \eqref{eq:jung3} and \eqref{eq:jung2} are not necessary.

\section{Main Results}

We first introduce the concept of generalized metric which will be used in our proofs.

\begin{definition}\label{def:genmetric}
	For a nonempty set $X$, a function $d:X\times X\rightarrow[0,\infty]$ is called a generalized metric on $X$ if and only if satisfies
	\begin{enumerate}
		\item[\bfseries M1] $d(x,y)=0$ if and only if $x=y$,
		\item[\bfseries M2] $d(x,y)=d(y,x)$ for all $x,y\in X$,
		\item[\bfseries M3] $d(x,z)\leq d(x,y)+d(y,z)$ for all $x,y,z\in X$.
	\end{enumerate}
\end{definition}

\noindent It should be remarked that the only difference of the generalized metric from the usual metric is that the range of the former is permitted to be an unbounded interval.

We will use the following fixed point result as main tool in our proofs, we refer to \cite{fixpoint} for the proof of this result.

\begin{theorem}\label{th:fixpoint}
	Let $(X,d)$ be a generalized complete metric space. Assume that $T:X\rightarrow X$ is a strictly contractive operator with the Lipschitz constant $\Lambda<1$. If there is a nonnegative integer $k$ such that $d\left(T^{k+1}x,T^kx\right)<\infty$ for some $x\in X$, then the following are true:
	\begin{enumerate}
		\item[(a)] The sequence $\left\{T^nx\right\}$ converges to a fixed point $x^*$ of $T$,
		\item[(b)] $x^*$ is the unique fixed point of $T$ in $$X^*=\left\{y\in X \; : \; d\left(T^kx,y\right)<\infty \right\},$$
		\item[(c)] If $y\in X^*$, then $$d\left(y,x^*\right)\leq\frac{1}{1-\Lambda}d\left(Ty,y\right).$$
	\end{enumerate}
\end{theorem}

In our proofs, we will use the metric defined in the following Lemma. In order to be able to use Theorem \ref{th:fixpoint}, we will need completeness of the space $X:=C\left([a,b],\mathbb{R}\right)$ which is given in the following result (see \cite{ben1}).

\begin{lemma}\label{lm:completeness}
	Let $X:=C\left([a,b],\mathbb{R}\right)$ and define the function $d:X\times X\rightarrow[0,\infty]$ with
	\begin{equation*}\label{eq:dmetric}
		d\left(f,g\right):=\inf\left\{ C\in\left[0,\infty \right] \;:\; \left| f(t)-g(t) \right| \leq C\Phi(t),\; t\in [a,b] 
		\right\}
	\end{equation*}
	where $\Phi:[a,b]\rightarrow\left(0,\infty\right)$ is a given continuous function. Then $\left(X,d\right)$ is a generalized complete metric space.
\end{lemma}

For the whole of this section, we define the interval $I:=[t_0, t_0+r]$ for real numbers $t_0$ and $r$ with $r>0$, further we define the set $X$ of all continuous functions defined on $I$, i.e. $X:=C(I,\mathbb{R})$.

The following theorem is the first main result of this paper, which shows that the conditions \eqref{eq:jung3} and \eqref{eq:jung2} is not necessary for solutions of Volterra integral equation \eqref{eq1} to have Hyers-Ulam-Rassias stability on bounded intervals.

\begin{theorem}\label{th:rassias}
	Suppose that the function $f:I\times\mathbb{R}\to\mathbb{R}$ is a continuous function satisfying the Lipschitz condition
	\begin{equation}\label{eq:th2-lip}
	\left|f(t,y_1)-f(t,y_2)\right|\leq L\left| y_1-y_2 \right|
	\end{equation}
	for all $t\in I$, all $y_1,y_2\in\mathbb{R}$ and some $L>0$. If a continuous function $y:I\to\mathbb{R}$ satisfies
	\begin{equation*}\label{eq:th2-ineq1}
		\left|y(t)-\int_{t_0}^{t}f(s,y(s))\d{s}\right|\leq\varphi(t)
	\end{equation*}
	for all $t\in I$, where $\varphi:I\to(0,\infty)$ is a nondecreasing continuous function, then there exists a unique solution $y_0:I\to\mathbb{R}$ of Volterra integral equation \eqref{eq1} satisfying
	\begin{equation}\label{eq:th2-ineq3}
	\left|y(t)-y_0(t)\right|\leq\frac{\varphi(t){\rm e}^{\eta r}}{1-L/\eta}
	\end{equation}
	for all $t\in I$, where $\eta$ is an arbitrary fixed real number with $\eta>L$.
\end{theorem}

\begin{proof}
	For any fixed $\eta\in\mathbb{R}$ with $\eta>L$, let us introduce the generalized metric on $X$ by 
	\begin{equation}\label{eq:def-metric}
	d(f,g):=\inf\{ C\in[0,\infty]\;:\; \left| f(t)-g(t) \right|{\rm e}^{-\eta(t-t_0)}\leq C\varphi(t),\;t\in I \}.
	\end{equation}
	According to Lemma \ref{lm:completeness}, $(X,d)$ is a generalized complete metric space. Now define the operator $\Theta:X\to X$ by
	\begin{equation*}
		\left(\Theta y\right)(t):=\int_{t_0}^{t}f(s,y(s))\d{s}
	\end{equation*}
	for all $t\in I$.

	Since $\Theta y$ is continuously differentiable on $I$, we remark that $\Theta y\in X$. For a given $g_0\in X$, since $f$ and $g_0$ are bounded on $I$ and $\min_{t\in I}\varphi(t)>0$, there exists a constant $C_0<\infty$ such that
	$$\left|(\Theta g_0)(t)-g_0(t)\right|{\rm e}^{-\eta(t-t_0)}=\left|\int_{t_0}^{t}f(s,g_0(s))\d{s}-g_0(t)\right|{\rm e}^{-\eta(t-t_0)}\leq C_0\varphi(t),$$
	for all $t\in I$, which implies $d(\Theta g_0,g_0)<\infty$. Furthermore, for any give $g\in X$, since $g$ and $g_0$ are bounded on $I$ and $\min_{t\in I}\varphi(t)>0$, there exists a constant $C_1<\infty$ such that $$\left| g_0(t)-g(t) \right|{\rm e}^{-\eta(t-t_0)}\leq C_1\varphi(t)$$ for all $t\in I$, i.e. $d(g_0,g)<\infty$ and hence $\{ g\in X\,:\,d(g_0,g)<\infty \}=X$.

	Now we will show that $\Theta:X\to X$ is strictly contractive on $X$. First observe that, with integration by parts and using monotonicity of $\varphi$,
	\begin{equation}\label{eq:th:intbyparts}
	\int_{t_0}^{t}\varphi(s){\rm e}^{\eta(s-t_0)}\d{s}\leq \frac{\varphi(t)}{\eta}{\rm e}^{\eta(t-t_0)}-\frac{\varphi(t)}{\eta}\int_{t_0}^{t}\varphi'(s){\rm e}^{\eta(s-t_0)}\d{s}\leq\frac{\varphi(t)}{\eta}{\rm e}^{\eta(t-t_0)}
	\end{equation}
	for all $t\in I$. For any $g_1,g_2\in X$, let $C_{g_1,g_2}\in[0,\infty]$ be an arbitrary constant such that $d(g_1,g_2)\leq C_{g_1,g_2}$, i.e. $\left|g_1(t)-g_2(t)\right|{\rm e}^{-\eta(t-t_0)}\leq C_{g_1,g_2}\varphi(t)$ for all $t\in I$. Then we have, by using \eqref{eq:th2-lip} and \eqref{eq:th:intbyparts},
	
	\begin{eqnarray*}
		\left| (\Theta g_1)(t)-(\Theta g_2)(t) \right| &=& \left| \int_{t_0}^{t}f(s,g_1(s))\d{s}-\int_{t_0}^{t}f(s,g_2(s))\d{s} \right|\\
		&\leq&\int_{t_0}^{t}\left|f(s,g_1(s))-f(s,g_2(s))\right|\d{s}\\
		&\leq&L\int_{t_0}^{t}\left|g_1(s)-g_2(s)\right|\d{s}\\
		&=&L\int_{t_0}^{t}\left|g_1(s)-g_2(s)\right|{\rm e}^{-\eta(s-t_0)}{\rm e}^{\eta(s-t_0)}\d{s}\\
		&\leq&LC_{g_1,g_2}\int_{t_0}^{t}\varphi(s){\rm e}^{\eta(s-t_0)}\d{s}\\
		&\leq&\frac{L}{\eta}C_{g_1,g_2}\varphi(t){\rm e}^{\eta(t-t_0)}
	\end{eqnarray*}
	for all $g_1,g_2\in X$ and for all $t\in I$. Hence, we have $$\left| (\Theta g_1)(t)-(\Theta g_2)(t) \right|{\rm e}^{-\eta(t-t_0)}\leq \frac{L}{\eta}C_{g_1,g_2}\varphi(t),$$ for all $g_1,g_2\in X$ and all $t\in I$, that is $$d(\Theta g_1,\Theta g_2)\leq\frac{L}{\eta}d(g_1,g_2)$$ for all $g_1,g_2\in X$. Therefore, $\Theta:X\to X$ is strictly contractive and the assumptions of Theorem \ref{th:fixpoint} are satisfied with $k=1$ and $X^*=X$. According the Theorem \ref{th:fixpoint}, there exists a unique solution $y_0:I\to\mathbb{R}$ of Volterra integral equation \eqref{eq1} satisfying
	$$ d(y,y_0)\leq\frac{1}{1-\Lambda}d(\Theta y, y) \leq \frac{\varphi(t)}{1-L/\eta},$$ which implies $$\left|y(t)-y_0(t)\right|{\rm e}^{-\eta(t-t_0)}\leq \frac{\varphi(t)}{1-L/\eta}$$ for all $t\in I$. Thus, we have $$\left|y(t)-y_0(t)\right|\leq \frac{\varphi(t)}{1-L/\eta}{\rm e}^{\eta(t-t_0)}\leq\frac{\varphi(t){\rm e}^{\eta r}}{1-L/\eta}$$ for all $t\in I$, which completes the proof.
\end{proof}

As a corollary of Theorem \ref{th:rassias}, we obtain Hyers-Ulam stability of solutions of Volterra integral equation \eqref{eq1} on finite closed intervals without any restrictions on Lipschitz constant.

\begin{corollary}\label{cor:ulam}
	Suppose that the function $f:I\times\mathbb{R}\to\mathbb{R}$ is a continuous function satisfying the Lipschitz condition \eqref{eq:th2-lip} for all $t\in I$, all $y_1,y_2\in\mathbb{R}$ and some $L>0$. If a continuous function $y:I\to\mathbb{R}$ satisfies
	\begin{equation*}
		\left|y(t)-\int_{t_0}^{t}f(s,y(s))\d{s}\right|\leq\varepsilon
	\end{equation*}
	for all $t\in I$, then there exists a unique solution $y_0:I\to\mathbb{R}$ of Volterra integral equation \eqref{eq1} satisfying
	\begin{equation}\label{eq:cor:1}
	\left|y(t)-y_0(t)\right|\leq \frac{\varepsilon{\rm e}^{\eta r}}{1-L/\eta}
	\end{equation}
	for all $t\in I$, where $\eta$ is an arbitrary fixed real number with $\eta>L$.
\end{corollary}

\begin{remark}
	In Theorem \ref{th:rassias} we can obtain Hyers-Ulam-Rassias stability of solutions of Volterra integral equation \eqref{eq1} without the assumptions  $$\left| \int_{a}^{t}\varphi(s)\d{s} \right|\leq K\varphi(t)\qquad\textrm{and}\qquad KL<1,$$ while they are required in the result given in \cite{jung2}. We also note that result given in Theorem 3.1 on Hyers-Ulam stability of the paper \cite{jung2} works only under the assumption $0<Lr<1$, while our result Corollary \ref{cor:ulam} is valid for all $L>0$.
\end{remark}

Now we will give our second main result, which concerns the Ulam type stability of \eqref{eq2}. In a similar way, in the following result, we show that the conditions \eqref{eq:jung3} and \eqref{eq:jung2} is not necessary for solutions of Volterra integral equation \eqref{eq2} to have Hyers-Ulam-Rassias stability on bounded intervals.

\begin{theorem}\label{th:rassias2}
	Suppose that the function $f:I\times I\times\mathbb{R}\to\mathbb{R}$ is a continuous function satisfying the Lipschitz condition
	\begin{equation}\label{eq:th3-lip}
	\left|f(t,y, z_1)-f(t,y, z_2)\right|\leq L\left| z_1-z_2 \right|
	\end{equation}
	for all $t\in I$, all $z_1,z_2\in\mathbb{R}$ and some $L>0$. If a continuous function $y:I\to\mathbb{R}$ satisfies
	\begin{equation*}\label{eq:th3-ineq1}
		\left|y(t)-\int_{t_0}^{t}f(t, s,y(s))\d{s}\right|\leq\varphi(t)
	\end{equation*}
	for all $t\in I$, where $\varphi:I\to(0,\infty)$ is a nondecreasing continuous function, then there exists a unique solution $y_0:I\to\mathbb{R}$ of Volterra integral equation \eqref{eq2} satisfying \eqref{eq:th2-ineq3} for all $t\in I$, where $\eta$ is an arbitrary fixed real number with $\eta>L$.
\end{theorem}

\begin{proof}
	Proof of this result is very similar to proof of Theorem \ref{th:rassias}, hence we omit the details and give only the sketch of the proof.
	
	For any fixed $\eta\in\mathbb{R}$ with $\eta>L$, define the generalized metric on $X$ by \eqref{eq:def-metric}. According to Lemma \ref{lm:completeness}, $(X,d)$ is a generalized complete metric space. Now define the operator $\Theta:X\to X$ by
	\begin{equation*}
		\left(\Theta y\right)(t):=\int_{t_0}^{t}f(t,s,y(s))\d{s}
	\end{equation*}
	for all $t\in I$. As in proof of Theorem \ref{th:rassias}, one can show that $d\left(\Theta g_0,g_0\right)<\infty$ for all $g_0\in X$ and $\{ g\in X\,:\,d(g_0,g)<\infty \}=X$. Now we will show that the operator $\Theta:X\to X$ is strictly contractive on $X$. For any $g_1,g_2\in X$, let $C_{g_1,g_2}\in[0,\infty]$ be a constant such that $d(g_1,g_2)\leq C_{g_1,g_2}$, i.e. $\left|g_1(t)-g_2(t)\right|$ $\times{\rm e}^{-\eta(t-t_0)} \leq C_{g_1,g_2}\varphi(t)$ for all $t\in I$. Then we have, by using \eqref{eq:th3-lip} and \eqref{eq:th:intbyparts},
	
	\begin{eqnarray*}
		\left| (\Theta g_1)(t)-(\Theta g_2)(t) \right| &=& \left| \int_{t_0}^{t}f(t,s,g_1(s))\d{s}-\int_{t_0}^{t}f(t,s,g_2(s))\d{s} \right|\\
		&\leq&\int_{t_0}^{t}\left|f(s,g_1(s))-f(s,g_2(s))\right|\d{s}\\
		&\leq&L\int_{t_0}^{t}\left|g_1(s)-g_2(s)\right|\d{s}\\
		&=&L\int_{t_0}^{t}\left|g_1(s)-g_2(s)\right|{\rm e}^{-\eta(s-t_0)}{\rm e}^{\eta(s-t_0)}\d{s}\\
		&\leq&LC_{g_1,g_2}\int_{t_0}^{t}\varphi(s){\rm e}^{\eta(s-t_0)}\d{s}\\
		&\leq&\frac{L}{\eta}C_{g_1,g_2}\varphi(t){\rm e}^{\eta(t-t_0)}
	\end{eqnarray*}
	for all $g_1,g_2\in X$ and for all $t\in I$. Hence, we have $$\left| (\Theta g_1)(t)-(\Theta g_2)(t) \right|{\rm e}^{-\eta(t-t_0)}\leq \frac{L}{\eta}C_{g_1,g_2}\varphi(t),$$ for all $g_1,g_2\in X$ and all $t\in I$, that is $$d(\Theta g_1,\Theta g_2)\leq\frac{L}{\eta}d(g_1,g_2)$$ for all $g_1,g_2\in X$. Therefore, $\Theta:X\to X$ is strictly contractive and the assumptions of Theorem \ref{th:fixpoint} are satisfied with $k=1$ and $X^*=X$.
	
	As in the proof of Theorem \ref{th:rassias}, by using Theorem \ref{th:fixpoint}, we conclude that there exists a unique solution $y_0:I\to\mathbb{R}$ of Volterra integral equation \eqref{eq2} satisfying
	$$\left|y(t)-y_0(t)\right|\leq\frac{\varphi(t){\rm e}^{\eta r}}{1-L/\eta}$$ for all $t\in I$.
\end{proof}

As a corollary of Theorem \ref{th:rassias2}, by taking $\varphi(t):=\varepsilon$, we obtain Hyers-Ulam stability of solutions of Volterra integral equation \eqref{eq2} on finite closed intervals without the restriction \eqref{eq:jung1}.

\begin{corollary}\label{cor:ulam2}
	Suppose that the function $f:I\times I\times\mathbb{R}\to\mathbb{R}$ is a continuous function satisfying the Lipschitz condition \eqref{eq:th3-lip} for all $t\in I$, all $z_1,z_2\in\mathbb{R}$ and some $L>0$. If a continuous function $y:I\to\mathbb{R}$ satisfies
	\begin{equation*}
		\left|y(t)-\int_{t_0}^{t}f(t,s,y(s))\d{s}\right|\leq\varepsilon
	\end{equation*}
	for all $t\in I$, then there exists a unique solution $y_0:I\to\mathbb{R}$ of Volterra integral equation \eqref{eq2} satisfying \eqref{eq:cor:1} for all $t\in I$, where $\eta$ is an arbitrary fixed real number with $\eta>L$.
\end{corollary}

\section{Examples}

\begin{example}\label{ex:1}
	Consider the Volterra integral equation
	\begin{equation}\label{eq:ex1}
	y(t)=\int_{0}^{t}sy(s)\d{s}
	\end{equation}
	on the interval $I:=[0,2]$. The function $f(t,y(t))=ty(t)$ satisfies the Lipschitz condition on $I$ with the Lipschitz constant $L=2$ since $$\left|f(t,y_1)-f(t,y_2)\right|=t\left|y_1-y_2\right|\leq2\left|y_1-y_2\right|$$ for all $t\in I$. According to Corollary \ref{cor:ulam}, the Volterra integral equation \eqref{eq:ex1} is stable in the sense of Hyers-Ulam on the interval $I$. We remark that Theorem 3.1 of \cite{jung2} and Theorem 5.1 of \cite{castro} do not work in this problem since $\lambda r=2\cdot1=2>1$.
\end{example}

\begin{example}
	Consider the Volterra integral equation \eqref{eq:ex1} of Example \ref{ex:1} on the interval $I:=[0,2]$, it has been shown that the function $f(t,y(t))$ satisfies the Lipschitz condition with $L=2$. If we choose $\varphi(t):={\rm e}^t$, we have $$\left|\int_{0}^{t}\varphi(s)\d{s}\right|=\int_{0}^{t}{\rm e}^s\d{s}={\rm e}^t-1\leq{\rm e}^t=\varphi(t)$$ for all $t\in I$, that is, the inequality \eqref{eq:jung3} holds with $K=1$. According to Theorem \ref{th:rassias}, the Volterra integral equation \eqref{eq1} is stable in the sense of Hyers-Ulam-Rassias. We remark that the condition \eqref{eq:jung3} is not required in our result but this example shows that  even if the inequality \eqref{eq:jung3} holds, Theorem 2.1 of the paper \cite{jung2} and Theorem 3.1 of the paper \cite{castro} do not work on this problem since $KL=2>1$ in this case.
\end{example}

\bibliography{refs}

\end{document}